\documentclass[twocolumn]{autart}    % Enable this line and disable the 
\usepackage{graphicx}          % Include this line if your 
\usepackage[T1]{fontenc}
\usepackage{amssymb,amsmath}
\usepackage{color}
\usepackage{enumerate}
\begin{document}

\newcommand\mcrit{10} %critical value of m, for which 1-cycle is unique
\newcommand\mcritnext{\the\numexpr\mcrit+1\relax} %critical value + 1
\newcommand\mcritprev{\the\numexpr\mcrit-1\relax} %critical value + 1

%for some reasons, does not work correctly with Overleaf Tex + PDF. I had to redefine it.
\renewcommand\triangleq{\overset{\scriptscriptstyle\Delta}{=}}
\renewcommand\square{\Box}

\begin{frontmatter}

\title{%The Existence of 
Cycles in Impulsive Goodwin's Oscillators of Arbitrary Order\thanksref{footnoteinfo}} % Title, preferably not more 

\thanks[footnoteinfo]{This paper was not presented or submitted to elsewhere. Corresponding author A.V.~Proskurnikov.}

\author[PdT]{Anton V. Proskurnikov}\ead{anton.p.1982@ieee.org},    % Add the 
\author[UU]{H{\aa}kan Runvik}\ead{hakan.runvik@it.uu.se},               % e-mail address 
\author[UU]{Alexander Medvedev}\ead{alexander.medvedev@it.uu.se}  % (ead) as shown

\address[PdT]{Department of Electronics and Telecommunications, Polytechnic University of Turin, 10129, Italy}  % Please supply                                 
\address[UU]{Information Technology, Uppsala University, SE-752 37, Sweden}             % full addresses
%\address[Baiae]{The White House, Baiae}        % here.

                             % keyword list or with the 
                                          % help of the Automatica 
                                          % keyword wizard

\begin{abstract}                          % Abstract of not more than 200 words.
Existence of periodical solutions, i.e. cycles, in the Impulsive Goodwin's Oscillator (IGO) with the continuous part of an arbitrary order $m$ is considered. The original  IGO with a third-order continuous part is a hybrid model that portrays a chemical or biochemical system composed of three substances represented by their concentrations and arranged in a cascade. The first substance in the chain is introduced via an impulsive feedback where both the impulse frequency and weights are modulated by the measured output of the continuous part. It is shown that, under the standard assumptions on the IGO, a positive periodic solution with one firing of the pulse-modulated feedback in the least period also exists  in models with any $m\geq 1$.
Furthermore,  the uniqueness of this 1-cycle is proved for the IGO with $m\leq \mcrit$ whereas, for $m>\mcrit$, the uniqueness can still be guaranteed under mild assumptions on the frequency modulation function.
\end{abstract}

\end{frontmatter}

\section{Introduction}
 Analyzing the dynamics of systems that simultaneously operate in fast and slow time scale (slow-fast systems) is a classical problem leading to the theory of singularly-perturbed dynamical systems \cite{S85}. Fast dynamics, i.e. rapid evolution occurring over shorter times, can be approximated by the impact of  finite or infinite impulse sequences resulting in (state vector) jumps \cite{KL61}. The impulsive action is then modeled either as a feedback or independent discrete process, i.e. a realization of a Markov chain. In the former case, one deals with a pulse-modulated feedback \cite{GC98} or  event-triggered control \cite{HJT12}, whereas the latter leads to hybrid control with Markovian switching. 

Theory of impulsive differential equations \cite{LBS} constitutes the mathematical ground of impulsive systems analysis and design. Impulsive models organically arise in biomedical, mechanical, ecological, environmental applications and are  present virtually in all fields of science where mathematical modeling is utilized. Predator-pray models with application to, e.g. pest control, make use of impulsive signals to represent human action \cite{ZL05}. Impulses (impacts) appear in non-smooth mechanics due to hard constraints on state variables and control signals. Numerous examples of practically important mechanical systems with impacts, including gear boxes, railway bogie, vibration table,  are provided in \cite{P00}.   Periodical medical pharmacological treatments is another significant application area of impulsive dynamical systems, where modeling is typically aimed at optimizing the treatment protocol, \cite{CCP20}. Impulses reflect the way drugs are administered, namely through injections or orally in tablet formulation. Pulsatile mode of drug administration also arises when a physiological behavior is mimicked by a treatment. A profound example of this concept is the pulsatile artificial pancreas. The physiological regulation exercised via the pancreas during a meal results in a series of insulin pulses whose frequency and amplitude are modulated by the blood glucose level \cite{BTH17}.  Therefore, there is increasing interest in impulsive control of the artificial pancreas \cite{HLS12}.

 Impulsive systems possess non-smooth dynamics and, thus, can exhibit complex nonlinear behaviors. 
Solutions converging to an equilibrium or an oscillative attractor are observed in linear time-invariant (LTI) systems under pulse-modulated feedback. The latter can correspond to sustained periodic or non-periodic (chaotic, quasiperiodic) solutions. The impulsive Goodwin's oscillator (IGO) \cite{MCS06}, \cite{Aut09} is a hybrid system that generalizes the classical continuous Goodwin's oscillator \cite{Good65} by substituting the original continuous static nonlinear feedback with a pulse-modulated one. The IGO lacks equilibria and admits solutions of high periodicity as well as chaotic and quasiperiodic ones \cite{ZCM12b}.

The rationale behind the IGO was originally to incorporate the experimentally observed principle of pulsatile endocrine regulation \cite{WTT10} into a widely used mathematical model of biochemical oscillation. However, the IGO can be seen as a general construct that illustrates how sustained oscillation can be obtained in a positive (continuous) LTI system by means of positive-valued feedback, no matter what the nature of the system is. From that perspective, the dynamics of the continuous part have to be as nonrestrictive as possible. Yet, in previous work on the IGO, only first-order  \cite{ZCM12a} and third-order continuous LTI dynamics have been addressed. In the latter case, the continuous dynamics augmented with point-wise \cite{CMM14} or distributed delay \cite{CM16} were also considered. 

The present paper generalizes the IGO structure to continuous LTI blocks of higher order than three. From an application point view, the order of an LTI model is a degree of freedom exploited by the designer to obtain a parsimonious description of essential model properties. Then setting the model order to a fixed constant is impractical. Further, when the model variables correspond to physical or chemical properties, the model order is defined by the number of variables whose time evolution has to be captured. Naturally, the number of dynamically interacting quantities in a concrete system can be arbitrary large.

Sustained rhythmical  behaviors are ubiquitous in nature \cite{GM88}. It is debatable whether such a behavior is suitably modeled as a perturbed  periodic solution of a dynamical system or a chaotic such. In the IGO, the main bifurcation mechanism  leading till chaos is frequency doubling \cite{ZCM12b}. Therefore, the existence of a periodic solution  is a central question in the IGO as it defines its very function. The focus here is, consequently, on a simplest kind of periodic solution that is characterized by just one impulse in the pulse-modulation feedback in the least period, i.e. a 1-cycle.

%\subsubsection*{Contribution of this paper.}

In this paper,  a generalization of the IGO to models with arbitrary continuous part order $m$, henceforth termed as IGO($m$), is proposed. The existence and uniqueness of a 1-cycle in IGO($3$) were established in \cite{Aut09}. 
Further, the same properties were also proved for IGO($1$)~\cite{ZCM12a}. Here, we generalize this result to IGO($m$) both regarding cycle existence and uniqueness, which constitutes the main contribution of this work.

First, we show that IGO($m$) possesses at least one 1-cycle (Theorem~\ref{thm.exist}). Furthermore, this property applies to a broad class of impulsive systems with Hurwitz stable and positive continuous-time part (Remark~\ref{rem.relax}).

Second, we prove that the 1-cycle is unique for dimensions $m\leq\mcrit$ (Theorem~\ref{thm.m9}), thus generalizing  Theorem~1 in~\cite{Aut09}.
As discussed in Section~\ref{sec.solution}, this development is far from being straightforward. It relies on the theory of divided differences and the Opitz formula allowing to compute an analytic function of a matrix with two-diagonal structure.

Third, we examine the problem of 1-cycle uniqueness in IGO($m$) with $m\geq\mcritnext$. Surprisingly, in this situation, the uniqueness may fail to hold and an example of such a case is given in Section~\ref{sec.num}. The uniqueness is, however, ensured if the derivative of the frequency modulator function does not attain anomalously large values (Theorem~\ref{thm.m10}). 

%\subsubsection*{The paper organization.}

The rest of the paper is organized as follows. After summarizing the notation, the IGO($m$) model is introduced in Section~\ref{sec.igo}.
Section~\ref{sec.setup} formulates the problem at hand, namely the existence and the uniqueness of 1-cycles in IGO($m$). Solutions to these problems are presented in Section~\ref{sec.solution}, with the proofs following separately in Section~\ref{sec.proof}. %We show that the uniqueness of 1-cycle always takes place when $m\leq\mcrit$ whereas, for $m>\mcrit$, an additional constraint on the frequency modulator function is needed, which generally cannot be discarded: 
An example of IGO($11$) with three distinct 1-cycles is given in Section~\ref{sec.num}. Appendices contains necessary information about divided differences and the Opitz formula (Appendix~\ref{sec.app-opitz}) and a proof of a technical lemma (Appendix~\ref{sec.app-psi}).

\section*{Notation}

The symbol $\triangleq$ henceforth means ``defined as''.

As usual, $\mathbb{R}$ and $\mathbb{R}_+$ stand, respectively, for the sets of all and nonnegative real numbers. The real vector space of dimension $m$ is then $\mathbb{R}^m$.  We use $\mathbb{N}_0$ to denote the set of nonnegative integers $\{0,1,\ldots\}$.

Given a function $f(\xi)$ of a scalar argument $\xi\in\mathbb{R}$, we denote its derivative evaluated at $\xi=\xi_0$ by $f'(\xi_0)$; $f^{(k)}(\xi)$ denotes the $k$th-order derivative. For a function of time $f(t)$, the derivative is equivalently denoted by $\dot f(t)$.  For a mapping $Q:\mathbb{R}^m\to\mathbb{R}^m$, the symbol $Q'(x)$ denotes the Jacobian matrix evaluated at $x\in\mathbb{R}^m$.

\section{The impulsive Goodwin's oscillator}\label{sec.igo}
Consider a continuous-time autonomous system
\begin{equation}\label{eq:LTI}
    \dot x(t)=Ax(t),\quad y(t)=Cx(t)
\end{equation}
with the state $x\in \mathbb{R}^m$, the output $y(t)\in\mathbb{R}$, and  the state-space matrices structured as
\begin{equation}
A=\begin{bmatrix}\label{eq:matrices}
-a_1 & 0 & \hdots &  & 0\\
g_1 & -a_2 & 0 & &\vdots\\
0 & g_2 & -a_3 &&\\
\vdots & & \ddots & \ddots &\\
0 &\hdots && g_{m-1} &-a_m
\end{bmatrix},
C=
\begin{bmatrix}
0\\
0\\
\vdots\\
1
\end{bmatrix}^{\top}.
\end{equation}
Assuming positive $a_i, i=1,\dots,m$ and $g_i, i=1,\dots,m-1$, the matrix $A$ is both Hurwitz and Metzler. 

Introduce an infinite sequence of time instants $t_n>0,n\in \mathbb{N}_0$ generated by the recursion 
\begin{equation}\label{eq:Phi}
    t_{n+1} =t_n+T_n, \quad 
T_n =\Phi(y(t_n)).
\end{equation}
%where $\Phi: \mathbb{R}\to\mathbb{R}$ is a \emph{continuous non-decreasing} function.
The state vector of system \eqref{eq:LTI} undergoes jumps at the times $t_n$ governed by 
\begin{equation}\label{eq:F}
\begin{gathered}
    x(t_n^+)=x(t_n^-)+\lambda_n B, \, \lambda_n=F(y(t_n)), \\ 
    B^{\top}=
\begin{bmatrix}
1
&0
&\hdots
&0
\end{bmatrix}.
\end{gathered}
\end{equation}
Here $\Phi:\mathbb{R}\to\mathbb{R}$ and $F:\mathbb{R}\to\mathbb{R}$ are known functions. In impulsive control systems~\cite{GC98}, they are usually referred to  as the frequency and amplitude modulation function, respectively. Interpreting the jumps as events, impulsive feedback~\eqref{eq:Phi},\eqref{eq:F} can be seen as a \emph{self-triggered}~\cite{HJT12} controller, because the output of the system at time $t_n$ uniquely determines the subsequent jump instant $t_{n+1}$.

With $m=3$, model~\eqref{eq:LTI}-\eqref{eq:F} is known as the impulsive Goodwin's oscillator (IGO)~\cite{CMZh16}. Below, a generalization of the IGO to an arbitrary order $m$ of the continuous part~\eqref{eq:LTI}, i.e. IGO($m$), is treated. %To enable explicit referencing to the order of \eqref{eq:LTI}, the generalized model~\eqref{eq:LTI}-\eqref{eq:F} is henceforth referred to as IGO($m$).

Notice that $\Phi,F$ are not generally required to be continuous  to guarantee a unique solution to hybrid system~\eqref{eq:LTI}-\eqref{eq:F}. Nevertheless, their continuity will be assumed to prove the existence of periodic solutions.
Following~\cite{Aut09}, we also assume that \begin{equation}                             \label{eq:bounds}
0<\Phi_1\le \Phi(y)\le\Phi_2,\; 0<F_1\le F(y)\le F_2\;\forall y\geq 0,
\end{equation}
where $\Phi_1$, $\Phi_2$, $F_1$, $F_2$ are positive constant numbers. This entails a number of important properties of the IGO that are proved similarly to the case of $m=3$~\cite{Aut09,ZCM12b}. Namely, IGO($m$) is a \emph{positive} system also for any order $m$, i.e., for positive initial conditions 
$\forall i: x_i(0)>0$, the solution remains positive $\forall i: x_i(t)>0$. Furthermore, a solution $x(t), t\in \lbrack 0, \infty)$ admits the following ultimate bounds
\begin{gather}\label{eq.bounded1}
V_i\leq\liminf_{t\to\infty}x_i(t)\leq \limsup_{t\to\infty}x_i(t)\leq H_i,\\
\begin{gathered}\label{eq.bounded2}
V_1=\frac{F_1}{e^{a_1\Phi_2}-1},\;\;H_1=\frac{F_2}{1-e^{-a_1\Phi_1}},\\
V_{i}=\frac{g_{i-1}}{a_i}V_{i-1},\;\;
H_i=\frac{g_{i-1}}{a_i}H_{i-1}, \quad\forall i=2,\ldots,m.
\end{gathered}
\end{gather}
In this paper, we focus on \emph{periodic} solutions such that $x(t+T)=x(t)$, for some $T>0$. For such a solution, $\liminf$ and $\limsup$ in~\eqref{eq.bounded1} can be omitted.

%\textcolor{red}{Alexander: please polish this paragraph!}
Motivated by application to feedback endocrine regulation, additional monotonicity restrictions were imposed  on the frequency and amplitude modulation functions of IGO in \cite{Aut09}.  It was in particular assumed that $\Phi$ is \emph{non-decreasing} and $F$ is \emph{non-increasing}. These assumptions are consistent with the experimentally observed behavior of the pulse-modulated feedback loop in testosterone (Te) regulation~\cite{KV98,WIC87}. A decrease in the concentration of Te  increases both the frequency and amplitude of the gonadotropin-releasing hormone pulses, which in turn stimulate the Te production.
In fact, as will be shown in this paper, the existence of periodic solutions does not require the monotonicity assumption. Moreover, we will prove that a certain periodic solution termed as 1-cycle always exists. At the same time, the monotonicity allows to prove, under certain conditions, the uniqueness of 1-cycle.

\section{Problem formulation: 1-cycle}\label{sec.setup}

%The IGO is primarily a model of sustained oscillations, cf.~\eqref{eq:bounds}, and periodic solutions represent a simplest form of these.
A fundamental property of IGO($3$) established in~\cite{Aut09} is that it always possesses a unique 
periodic solution featuring only one jump over the (minimal) period $T>0$, i.e. a \emph{1-cycle}. Then, \eqref{eq:Phi} becomes
\[
t_{n+1}=t_n+T, \quad \Phi(y(t_n))=T, \quad \forall n\in \mathbb{N}_0.
\]
%which means, in view of~~\eqref{eq:Phi}, that $\Phi(y(t_n))=T$ for all $n$.
With the notation $X_n=x(t_n^-)$, the return map $X_{n+1}=Q(X_n)$, for $n=0,1,\ldots$,  of IGO($m$) is given \cite{Aut09,churilov2020} by
\begin{equation}\label{eq:Q}
Q(x)=e^{\Phi(Cx)A}(x+F(Cx)B),\quad x\in\mathbb{R}_{+}^m.
\end{equation}

%\paragraph*{1-cycles as fixed points}

As shown in~\cite{Aut09,churilov2020}, a 1-cycle corresponds to a \emph{fixed point} of the map $Q$. For such a point $x_*=Q(x_*)$, the corresponding 1-cycle is found as
\begin{equation}\label{eq.1-cycle-full}
\begin{gathered}
x(t)=e^{(t-t_n)A}(x_*+F(Cx_*)B),\quad t\in (t_n,t_{n+1}),\\
X_n=x_*,\; x(t_n^+)=x_*+F(Cx_*)B,\\
t_{n+1}=nT,\quad T=\Phi(Cx_*),\quad n\in \mathbb{N}_0.
\end{gathered}
\end{equation}

In view of the positivity of the IGO($m$), admissible 1-cycles correspond to fixed points $x_*\in\mathbb{R}_+^m$; For such a solution, periodic solution~\eqref{eq.1-cycle-full} will stay in $\mathbb{R}_+^m$.

In this paper, we address the problems of existence and uniqueness of feasible (positive) fixed points:\\
\textbf{Problem A.} \textit{Does IGO($m$) always have a feasible 1-cycle? Equvalently, does the corresponding mapping $Q$ have a fixed point $x_*=Q(x_*)\in\mathbb{R}_+^m$?}

%In~\cite{Aut09}, the existence was proven in the $m=3$ case. 
Below, in Theorem~\ref{thm.exist}, we give an affirmative answer to Problem~A for an \emph{arbitrary} $m$. This existence property is actually valid for a much more general class of impulsive systems (Remark~\ref{rem.relax}). 

A natural question of how many distinct 1-cycles an IGO($m$) might have then arises:\\
\textbf{Problem B.} \textit{Is the feasible 1-cycle of IGO($m$) (equivalently, the fixed point $x_*\in\mathbb{R}_+^m$ of $Q$) unique?}

The uniqueness of 1-cycle  for IGO($3$)  established in~\cite[Theorem~1]{Aut09} is generalized to $m\leq\mcrit$ in Theorem~\ref{thm.m9} of the present paper. For $m=\mcritnext$, however, it is possible to find  parameter values $a_i,g_i>0$ and functions $F,\Phi$ such that the corresponding IGO has three distinct 1-cycles. Multiple 1-cycles are although highly uncommon. As Theorem~\ref{thm.m9} and Theorem~\ref{thm.m10} show, the uniqueness can always be secured by limiting  $\Phi'$ or by letting the impulses to be sufficiently sparse, i.e. bounding $\Phi_1$ from below.

\section{Main Results}\label{sec.solution}

In this section, we state the main result of the paper providing complete solutions to Problem~A and Problem~B and formulated in Theorem~\ref{thm.exist}--Theorem~\ref{thm.m10}. Their proofs are summarized separately in Section~\ref{sec.proof}.
\subsection*{Problem A: Existence of 1-cycles in IGO($m$)}
%\paragraph*{Equation of periods and 1-cycles}
%We start with introducing 
The so-called ``equation of periods''~\cite{Aut09}  characterizes the feasible fixed points of $Q$ introduced in \eqref{eq:Q}
\begin{equation}\label{eq.period}
y=R(y)\triangleq F(y)C(e^{-\Phi(y)A}-I)^{-1}B,\quad y\in\mathbb{R}_+.
\end{equation}
Since $\Phi(y)>0$ for $y\geq 0$ thanks to~\eqref{eq:bounds} and $A$ is Hurwitz, the inverse matrix in~\eqref{eq.period} is well-defined.

Note that for $x\in\mathbb{R}^m_+$, the equation $Q(x)=x$ can be equivalently written as
\begin{equation}\label{eq.period0}
x=F(Cx)(e^{-\Phi(Cx)A}-I)^{-1}B,
\end{equation}
and, therefore, $y=Cx$ obeys~\eqref{eq.period}. 
Conversely, if $y$ is a root of equation~\eqref{eq.period}, then $x=F(y)(e^{-\Phi(y)A}-I)^{-1}B$ obeys~\eqref{eq.period0}, entailing  $Q(x)=x$. 
However, it is not obvious that such a vector $x$ is positive and the latter fact is ensured by one of the statements in Theorem~\ref{thm.exist} below.

%\subsection*{Problem A: Existence of 1-cycles in IGO($m$).}

%We now formulate our first main result.
\begin{thm}\label{thm.exist}
For all values of the parameters $a_1,\ldots,a_m>0$, $g_1,\ldots,g_{m-1}>0$ and \emph{continuous} functions $\Phi,F$ obeying~\eqref{eq:bounds}, the following statements are valid:
\begin{enumerate}
\item The function $R(\cdot)$ defined in~\eqref{eq.period} is uniformly strictly positive and bounded on $[0,\infty)$;
\item Equation~\eqref{eq.period} has at least one solution; all its solutions are strictly positive for $y>0$;
\item For every solution of~\eqref{eq.period}, the vector $x=F(y)(e^{-\Phi(y)A}-I)^{-1}B$ is a positive fixed point of  return map~\eqref{eq:Q};
\end{enumerate}
Hence, IGO($m$) always has at least one positive 1-cycle.
\end{thm}

Noticeably, Theorem~\ref{thm.exist} does not impose any monotonicity restrictions on $F$ and $\Phi$. 
As will be shown (Remark~\ref{rem.relax}), this theorem generalizes to a broad class of impulsive systems with positive and stable continuous-time part~\eqref{eq:LTI}, whose matrices $A,B,C$ may differ in structure from~\eqref{eq:matrices}. 

\subsection*{Problem B: Uniqueness of 1-cycles in IGO($m$)}

An elegant result established in~\cite[Theorem~1]{Aut09} states that, in the case $m=3$, the solution to~\eqref{eq.period} is \emph{unique}, because the function $R$ is \emph{non-increasing} on $[0,\infty)$.
This monotonicity property, proved in~\cite{Aut09} for $m=3$ by evoking the Jenssen inequality, remains valid for $1\leq m\leq\mcrit$, as shown below.
\begin{thm}\label{thm.m9}
For all $1\leq m\leq \mcrit$, positive parameter values $a_i,g_i>0$, and continuous non-increasing functions $F$ and non-decreasing functions $\Phi$ satisfying~\eqref{eq:bounds}, the function $R$ defined in~\eqref{eq.period} is
non-increasing on $[0,\infty)$. In particular,~\eqref{eq.period} has a unique positive solution, and the corresponding IGO($m$) has a unique 1-cycle. 

These statements retain their validity if one replaces the condition $m\leq\mcrit$ by the inequality
\begin{equation}\label{eq.phi1}
\frac{m-1}{\min_ia_i}\leq\Phi_1.
\end{equation}
\end{thm}

In Section~\ref{sec.proof}, we will show that the uniqueness of 1-cycle cannot be  generally established for $m=\mcritnext$ and it is possible to find an IGO($\mcritnext$) with at least three different 1-cycles. The numerical example of this in Section~\ref{sec.num} requires the function $\Phi$ to possess very large derivative at some points (violating also~\eqref{eq.phi1}). By forbidding excessive values of $\Phi'$, one can guarantee uniqueness for the order $m\geq\mcritnext$ as stated by the next theorem.

Recall the Riemann $\boldsymbol{\zeta}$-function
\begin{equation}\label{eq.zeta}
\boldsymbol{\zeta}(s)\triangleq\sum\nolimits_{k=1}^{\infty}k^{-s},\quad s>1.
\end{equation}
\begin{thm}\label{thm.m10}
Consider an IGO($m$), $m\geq\mcritnext$, whose modulation functions $\Phi,F$ are, respectively, non-increasing and non-decreasing. Assume also that $\Phi,F$ are absolutely continuous and, furthermore, for each $y>0$, one has
\begin{equation}\label{eq.phi-d-bound}
\Phi'(y)\leq\frac{C_m}{g_1\ldots g_{m-1}F(0)},\; C_m\triangleq\frac{(2\pi)^m}{2(m-1)\boldsymbol{\zeta}(m)}.
\end{equation}
Then, the function $y-R(y)$ is strictly increasing on $\mathbb{R}_+$, equation~\eqref{eq.period} has only one solution, 
and thus the IGO($m$) has a unique 1-cycle.
\end{thm}

%In view of Theorem~\ref{thm.m9}, condition \eqref{eq.phi-d-bound} in Theorem~\ref{thm.m10} is of interest for $m\geq\mcritnext$. 
The sequence $\boldsymbol{\zeta}(m), m=1,2,\dots$ is decreasing  and  $\boldsymbol{\zeta}(m)\geq 1$. Hence, for all $m\geq\mcritnext$, $\boldsymbol{\zeta}(m)\leq \boldsymbol{\zeta}(\mcritnext)\approx 1.005$ and, consequently, $C_m$ grows exponentially as $m\to\infty$. Numerical evaluation in Matlab yields $C_{11}\approx 3.01\cdot 10^7$, $C_{12}\approx 1.72\cdot 10^8$, 
$C_{13}\approx 9.91\cdot 10^8$. Condition~\eqref{eq.phi-d-bound} is thus not very restrictive for large $m$ but, yet,  cannot be fully abolished (see the example in Section~\ref{sec.num}). 

%\textcolor{blue}{I recall that $m=1$ was a pathological case. This should be suitably commented on. Here $C_m$ is infinite for $m=1$.}

\subsection*{Discussion}

%\paragraph*{The derivatives of $\Phi,F$ and stability of 1-cycles.}

Remarkably, none of Theorem~\ref{thm.exist}--Theorem~\ref{thm.m10} requires the nonlinearities $F,\Phi$ to be differentiable everywhere. In the case of Theorem~\ref{thm.m10}, we only need absolute continuity, which ensures existence of $\Phi'(y),F'(y)$ at \emph{almost every} point $y\in\mathbb{R}_+$). 
If $F$ and $\Phi$ are \emph{continuously differentiable} in a vicinity of the fixed point $y_*=R(y_*)$ in \eqref{eq.period}, then the (local exponential) \emph{orbital stability} of the corresponding 1-cycle can be examined (see~\cite{churilov2020}, where the underlying stability definitions can be found). Namely, 
the 1-cycle defined by a fixed point of the map $Q$ (as stated in Theorem~\ref{thm.exist}) %$x_0=F(y_0)(e^{-\Phi(y_0)A}-I)^{-1}B$, 
 is orbitally stable if and only if the Jacobian matrix $Q'(x_*)$ is Schur stable~\cite[Theorem~3]{churilov2020}. Obviously, $Q'(x_*)$ is fully determined by the parameters of continuous part~\eqref{eq:matrices} and the values $F(y_*),\Phi(y_*),F'(y_*),\Phi'(y_*)$.

%\paragraph*{The case of $m=3$~\cite{Aut09} versus our results.}

It should be noticed that the method of proving Theorem~1 in~\cite{Aut09}, although it yields the results of Theorem~\ref{thm.exist} and Theorem~\ref{thm.m9} for $m=3$, is not applicable to a general IGO($m$) for several reasons.
First, both existence and uniqueness are derived in~\cite{Aut09} from the monotonicity of the function $R$, which, as proved above, does not hold for $m>\mcrit$ without additional assumptions, whereas Theorem~\ref{thm.exist} (existence of 1-cycle) retains its validity. Second, the method of proving this monotonicity property is based on an analytic representation of $R(y)$ and its derivative $R'(y)$ in $m=3$ case, which was obtained in the proof of~\cite[Theorem~1]{Aut09} by a straightforward computation\footnote{In fact,~\cite{Aut09} adopts a modeling assumption that $a_1,a_2,a_3$ are pairwise distinct, which is abolished here.}. In the general case considered in the present paper, these two functions are computed by using the Opitz formula and the method of divided differences (Appendix~\ref{sec.app-opitz}), which tools are not exploited in~\cite{Aut09}.
Third, the closed-form representation of the derivative $R'(y)$ allows to derive its positivity from the Jenssen inequality, which trick, to the best of our knowledge, cannot be applied for $m>3$.
Hence, while following the same line of reasoning as~\cite{Aut09}, this paper substantially generalizes the  results of the latter by applying a different set of mathematical tools.

\section{Proofs and Auxiliary Results}\label{sec.proof}

This section summarizes the proofs of the Theorems formulated in Section~\ref{sec.solution} and also establishes  necessary auxiliary technical statements that might be of use elsewhere. %, i.e. Lemmas~\ref{lem.positivity},~\ref{lem.r-as-dd},~\ref{lem.Psi-props} and Corollaries~\ref{cor.positivity},~\ref{cor.monotone},~\ref{cor.r-d-est}.

\subsection{Lemmas and proof of Theorem~\ref{thm.exist}}

Recall that the matrix is called \emph{nonnegative} (respectively, Metzler) if all its entries (respectively, all its off-diagonal entries) are nonnegative.
Hence, if $M$ is a Metzler matrix, then $M+nI$ is nonnegative for some $n\in\mathbb{R}$ being large enough.

Following~\cite{HornJohnsonBook}, we introduce the \emph{graph} $\Gamma(B)$ of a nonnegative square matrix $B=(b_{ij})_{i,j\in I}$. In this graph, the nodes are in one-to-one correspondence with the elements of the index set $I$ and a directed arc $(i,j)$ is present if and only if $b_{ij}>0$. Positive diagonal entries stand for self-arcs. For each $k=1,2,\ldots$, the matrix $B^k$ has a positive entry $(B^k)_{ij}>0$ if and only if $\Gamma(B)$ contains a directed walk of length $k$ connecting $i$ to $j$ (this walk may contain self-loops and visit some vertices multiple times). 
We may formally generalize the definition of the graph to Metzler square matrices: given such a matrix $A=(a_{ij})_{i,j\in I}$, we connect two nodes $i,j\in I$ if and only if $a_{ij}>0$. 

The proof of Theorem~\ref{thm.exist} is based on the following positivity result.
\begin{lem}\label{lem.positivity}
For every Metlzer and Hurwitz matrix $A\in \mathbb{R}^{n\times n}$, the matrix-valued function $\Theta(\xi)=(e^{-\xi A}-I)^{-1}$ exists and is non-negative for all $\xi>0$ with strictly positive diagonal entries $\Theta_{ii}(\xi)>0\,\ i=1,\dots, n$.
Each off-diagonal entry $\Theta_{ij}(\xi),\,i\ne j$ is positive $\forall \xi>0$ if and only if
the graph $\Gamma[A]$ contains a path from $i$ to $j$; otherwise, $\Theta_{ij}(\xi)\equiv 0$.
\end{lem}
\begin{pf}
Notice first that the matrix exponential
\[
e^B=\sum_{k=0}^{\infty}\frac{1}{k!}B^k
\]
of a \emph{nonnegative} matrix $B$ is also a nonnegative matrix. Furthermore, $(e^B)_{ij}>0$ if and only if either $i=j$ or a directed walk from $i$ to $j$ exists in $\Gamma(B)$.

For an arbitrary $\xi\in\mathbb{R}$, the graphs of $A$ and $A+\xi I$ may differ only by the presence of self-arcs, which do not influence connectivity. Hence, two nodes $i$ and $j\ne i$ are connected by a directed walk in $\Gamma(A)$ if and only if they are connected by such in $\Gamma(A+\xi I)$. Choosing $\xi>0$ large enough, the matrix $A+\xi I$ is nonnegative. Hence, $(e^A)_{ij}=e^{-\xi}(e^{A+\xi I})_{ij}$ is nonnegative for all $i,j$, being positive if and only if $i=j$ or a walk leads from $i$ to $j$ in $\Gamma(A)$. 
The same statements hold true for $e^{\xi A}$ if $\xi>0$, because matrix $\xi A$ is Hurwitz and Metzler for any $\xi>0$, having same graph as $A$.
 
By noticing that $e^{\xi A}$ has the eigenvalues $e^{\xi\lambda_j(A)}$, where $\lambda_j(A)$ are the eigenvalues of $A$, one concludes that the exponential $e^{\xi A}$ has the spectral radius $\rho(e^{\xi A})=\max_j|e^{\lambda_j(A)}|=e^{\max_j\textrm{Re\,}\lambda_j(A)}<1$. Hence
\begin{equation*}%\label{eq.aux}
\Theta(\xi)=e^{\xi A}(I-e^{\xi A})^{-1}=e^{\xi A}\sum_{k=0}^{\infty}(e^{\xi A})^k=\sum_{k=1}^{\infty}(e^{\xi A})^{k}
\end{equation*}
is well-defined and nonnegative for all $\xi>0$. Furthermore, $\Theta_{ij}(\xi)>0$ if and only if 
$i=j$ or $i$ is connected to $j$ by a walk in $\Gamma(A)$ (otherwise, the $(i,j)$-entry of all summands vanishes), which completes the proof $\square$
\end{pf}
\begin{cor}\label{cor.positivity}
Let $A=(a_{ij})_{i,j\in I}$ be a Hurwitz and Metzler matrix and $b,c$ be two nonnegative column vectors of same dimension as $A$. Then 
$c^{\top}\Theta(\xi)b>0$ for all $\xi>0$ if and only if there exist indices $i,j\in I$ such that the elements $c_i>0,b_j>0$ and either $i=j$ or $\Gamma(A)$ contains
a directed walk from $i$ to $j$. If the latter condition is violated, then $c^{\top}\Theta(\xi)b\equiv 0$.
\end{cor}
\begin{pf}
Notice that $c^{\top}\Theta(\xi)b=\sum_{i,j}c_i\Theta_{ij}(\xi)b_j$. For $\xi>0$,
all summands in the latter sum are nonnegative, and thus $c^{\top}\Theta(\xi)b\geq 0$. The latter inequality is strict if and only if
at least one summand is positive $c_i\Theta_{ij}(\xi)b_j>0$, which is possible if and only if $c_i,b_j>0$ and $\Theta_{ij}(\xi)>0$. 
The statement now follows from Lemma~\ref{lem.positivity} $\square$
\end{pf}
\begin{cor}\label{cor.positivity1}
For the matrix $A$ in~\eqref{eq:matrices} and the column $B$ in \eqref{eq:F}, the column $(e^{-\xi A}-I)^{-1}B$ is positive for $\xi>0$.
\end{cor}
\begin{pf}
The graph $\Gamma(A)$ contains a unidirectional chain $n\to (n-1)\to\ldots\to 1$ thanks to inequalities $g_i>0\,\forall i$. Hence, each node $i=2,\ldots,n$
is connected to $1$ by a directed walk. Applying Corollary~\ref{cor.positivity} to $A$, $b=B$ and the coordinate vectors $c=e_1,e_2,\ldots,e_n$, one concludes that
all elements of $(e^{-\xi A}-I)^{-1}B$ are positive for $\xi>0$ $\square$
\end{pf}

\subsubsection*{Proof of Theorem~\ref{thm.exist}}

Now all the auxiliary results are in place to prove the claim of Theorem~\ref{thm.exist}.

Corollary~\ref{cor.positivity1} ensures that function $r(\xi)=C(e^{-\xi A}-I)^{-1}B$ is positive for $\xi>0$. Also, $r(\xi)$ is continuous at every point $\xi\in (0,\infty)$.
Notice that $R(y)=F(y)r(\Phi(y))$ due to~\eqref{eq.period}. In view of~\eqref{eq:bounds}, for every $y\geq 0$, one has
\[
0<F_1\min_{\xi\in[\Phi_1,\Phi_2]}r(\xi)\leq R(y)\leq F_2\max_{\xi\in[\Phi_1,\Phi_2]}r(\xi)<\infty,
\]
where the minimum and the maximum exist due to the Weierstrass extreme value theorem. This proves Statement~(1) of the Theorem.

Recalling that $F,\Phi$ are assumed to be continuous, $y-R(y)$ is a continuous function on $\mathbb{R}_+$ attaining a negative value at $y=0$
and positive values where $y$ is large enough. Statement (2) is now straightforward from the intermediate value theorem.
Statement~(3) follows from Corollary~\ref{cor.positivity1}, recalling that $F(y)\geq F_1>0$ $\square$

%\subsubsection*{Additional remarks}

\begin{rem}\label{rem.bounds}
{\rm
The proof of Theorem~\ref{thm.exist} implies that all the roots of~\eqref{eq.period} belong, in fact, to the closed interval
\[
F_1\min_{\xi\in[\Phi_1,\Phi_2]}r(\xi)\leq y=R(y)\leq F_2\max_{\xi\in[\Phi_1,\Phi_2]}r(\xi).
\]
The minimum and maximum can, in turn, be estimated by using the explicit representation of $r$ provided by \eqref{eq.r} and Lemma~\ref{lem.mean.val} in Section~\ref{thm.unique}.
This facilitates the numerical solution of~\eqref{eq.period} by e.g. the bisection method. We omit the technical details here for brevity.
}
\end{rem}
\begin{rem}\label{rem.relax}
{\rm
Theorem~\ref{thm.exist} can be generalized to guarantee the existence of a 1-cycle in a broader class of impulsive systems~\eqref{eq:LTI},~\eqref{eq:Phi},~\eqref{eq:F} than those with the matrix structures  specified in \eqref{eq:matrices}.
Corollary~\ref{cor.positivity1} remains valid for any Hurwitz and Metzler matrix $A$ and a column $B$ such
that each node $i$ of the graph $\Gamma(A)$ either corresponds to $B_i>0$ or is connected by a path to some node $j$ such that $B_j>0$. In such a case, statements (1)-(3) remain valid.
}
\end{rem}
%One may also notice that the monotonicity properties of $F,\Phi$ are unimportant. Statement~(1) relies only on %bounds~\eqref{eq:bounds}, whereas statements~(2) and~(3) require additionally that $F,\Phi$ are continuous. \textcolor{blue}%{Irrelevant to the existence of 1-cycle but impact other properties. Monotonicity is important for positive-valued negative feedback %and stability of fixed points.}
\begin{rem}
{\rm
Notice also that, assuming that $F,\Phi$ are continuous and~\eqref{eq:bounds} holds, the map $Q$ admits a nonnegative fixed point (whose components, however, may be zero) whenever 
$A$ is Hurwitz and Metzler and $B,C$ are nonnegative. Indeed, Corollary~\ref{cor.positivity} states that either $r(\xi)>0\,\forall\xi>0$ or $r(\xi)\equiv 0$ (in this degenerate case, the output $y(t)$ is decoupled from the input $u(t)$).
In the former case, Statement~(1) and Statement~(2) of Theorem~\ref{thm.exist} are valid; in the latter case, the equation of periods $y=R(y)$
has the unique solution $y=0$. In both cases, vector $x=F(y)(e^{-\Phi(y) A}-I)^{-1}B$ is a nonnegative 
(Lemma~\ref{lem.positivity}) fixed point of the map $Q$.
}
\end{rem}

\subsection{Lemmas and proofs of Theorem~\ref{thm.m9}--Theorem~\ref{thm.m10}}\label{thm.unique}

In this subsection, we intensively use divided differences (DD) and the Opitz formula (see Appendix~\ref{sec.app-opitz}  where the necessary background is summarized). 

Given a function $f:I\to\mathbb{R}$ on the interval $I\subseteq\mathbb{R}$ and $k+1$ points $x_0,\ldots,x_k$,  $k\in \mathbb{N}_0$,
$f[x_0,\ldots,x_k]$ stands for the $k$-th order DD (briefly, $k$-DD) evaluated at $x_0,\ldots,x_k$.

A useful property of $k$-DD in the present context is the following extension of the mean-value theorem.
\begin{lem}~\cite[Section~8]{DeBoor2005}\label{lem.mean.val}
Suppose that $f:I\to\mathbb{R}$ is $k$ times differentiable on $I$ and let $x_0,\ldots,x_k\in I$.
Then a point $\bar x\in[\min_i x_i,\max_i x_i]$ exists such that
\begin{equation}\label{eq.mean-value}
f[x_0,\ldots,x_k]=\frac{1}{k!}f^{(k)}(\bar x).
\end{equation}
By substituting $x_0=\ldots=x_k=\xi$, one thus has
\begin{equation}\label{eq.derivative}
f[\underbrace{\xi,\ldots,\xi}_{k+1}]=\frac{1}{k!}f^{(k)}(\xi)\;\forall\xi\in I.
\end{equation}
\end{lem}

%\subsubsection*{Functions $r(\xi)$ and $r'(\xi)$ as DDs}
Consider the function
$$r(\xi)\triangleq C(e^{-\xi A}-I)^{-1}B.$$
Lemma~\ref{lem.r-as-dd} provides representations of $r(\xi)$ and its derivative in terms of DD.  These results are instrumental in proving Theorem~\ref{thm.m9} and Theorem~\ref{thm.m10}.

Introduce two auxiliary functions
\[
\varphi(x)=1/(e^{x}-1),\quad \psi(x)=-x\varphi'(x)=xe^{x}/(e^{x}-1)^2.
\]
The Opitz formula (Appendix~\ref{sec.app-opitz}, Equation~\ref{eq.opitz}) leads to the following  result.
\begin{lem}\label{lem.r-as-dd}
Consider state-space matrices \eqref{eq:matrices}, and let $\bar g=g_1\ldots g_{m-1}>0$. Then, 
for all $\xi>0$, the function $r(\xi)$ and its derivative are found as
\begin{gather}
r(\xi)=(-\xi)^{m-1}\bar g\varphi[\xi a_1,\ldots,\xi a_m],\label{eq.r}\\
r'(\xi)=(-\xi)^{m-2}\bar g\psi[\xi a_1,\ldots,\xi a_m]\label{eq.r-d}.
\end{gather}
\end{lem}
\begin{pf}
\begin{description}
\item[Step 1:] Let $\bar g_i\triangleq g_1\ldots g_{i-1}$ for $i=2,\ldots,m$, (hence, $\bar g_m=\bar g$), and $\bar g_1\triangleq 1$.
Notice first that
\[
A=S\Lambda S^{-1},\quad S={\rm diag}(\bar g_1,\bar g_2,\bar g_3,\ldots,\bar g_{m}).
\]
The matrix $\Lambda$ is two-diagonal with the eigenvalues $-a_i, i=1, \dots, m$ on the main diagonal and ones on
the diagonal below. One can check that
\[
CS=\bar g_mC=\bar gC,\quad S^{-1}B=\bar g_1B=B.
\]
\item[Step 2:] For a function $f$ analytic in a vicinity of the eigenvalues $-a_1,\ldots,-a_m$,
one thus has $f(A)=Sf(\Lambda)S^{-1}$, furthermore, $Cf(A)B=\bar g f(\Lambda)_{m,1}$ (the subscript denotes the $(m,1)$ entry of the matrix $f(\Lambda)$). 
In virtue of~\eqref{eq.opitz}, it follows
\[
Cf(A)B=\bar gf[-a_1,\ldots,-a_m].
\]
Then Lemma~\ref{lem.scale} implies
\[
Cf(-\xi A)B=Cf_{-\xi}(A)B=(-\xi)^{m-1}\bar gf[\xi a_1,\ldots,\xi a_m],
\]
(in accordance with Lemma~\ref{lem.scale}, $f_{-\xi}(x)=f(-\xi x)$).
\item[Step 3:] Equality~\eqref{eq.r} is now straightforward by noticing that $r(\xi)=C\varphi(-\xi A)B$.
To prove~\eqref{eq.r-d}, recall that, for any differentiable invertible matrix function $X$, one has
$(X(\xi)^{-1})'=-X(\xi)^{-1}X'(\xi)X(\xi)^{-1}$. Therefore,
\[
\begin{split}
\frac{d}{d\xi}(e^{-\xi A}-I)^{-1}=(e^{-\xi A}-I)^{-1}Ae^{-\xi A}(e^{-\xi A}-I)^{-1}=\\
=Ae^{-\xi A}(e^{-\xi A}-I)^{-2}=(-\xi)^{-1}\psi(-\xi A).
\end{split}
\]
Hence, $r'(\xi)=(-\xi)^{-1}C\psi(-\xi A)B$, entailing~\eqref{eq.r-d} $\square$
\end{description}
%\textbf{Step 1.} 
%For each $i=2,\ldots,m$, let $\bar g_i\triangleq g_1\ldots g_{i-1}$ (hence, $\bar g_m=\bar g$), and $\bar g_1\triangleq 1$.
%Notice first that
%\[
%A=S\Lambda S^{-1},\quad S={\rm diag}(\bar g_1,\bar g_2,\bar g_3,\ldots,\bar g_{m}).
%\]
%Matrix $\Lambda$ is two-diagonal with the eigenvalues $-a_i$ on the main diagonal and ones on
%the diagonal below. One can check that
%\[
%CS=\bar g_mC=\bar gC,\quad S^{-1}B=\bar g_1B=B.
%\]
%\textbf{Step 2.} For each function $f$ analytic in a vicinity of the eigenvalues $-a_1,\ldots,-a_m$,
%one thus has $f(A)=Sf(\Lambda)S^{-1}$, furthermore, $Cf(A)B=\bar g f(\Lambda)_{n1}$.
%Using~\eqref{eq.opitz}, one shows that
%\[
%Cf(A)B=\bar gf[-a_1,\ldots,-a_m].
%\]
%Using Lemma~\ref{lem.scale}, one shows that
%\[
%Cf(-\xi A)B=Cf_{-\xi}(A)B=(-\xi)^{m-1}\bar gf[\xi a_1,\ldots,\xi a_m].
%\]
%Here, in accordance with Lemma~\ref{lem.scale}, $f_{-\xi}(x)=f(-\xi x)$.
%\textbf{Step 3.} Equation~\eqref{eq.r} is now straightforward by noticing that $r(\xi)=C\varphi(-\xi A)B$.
%To prove~\eqref{eq.r-d}, recall that, for any differentiable invertible matrix function $X$, one has
%$(X(\xi)^{-1})'=-X(\xi)^{-1}X'(\xi)X(\xi)^{-1}$. Therefore,
%\[
%\begin{split}
%\frac{d}{d\xi}(e^{-\xi A}-I)^{-1}=(e^{-\xi A}-I)^{-1}Ae^{-\xi A}(e^{-\xi A}-I)^{-1}=\\
%=Ae^{-\xi A}(e^{-\xi A}-I)^{-2}=(-\xi)^{-1}\psi(-\xi A).
%\end{split}
%\]
%Hence, $r'(\xi)=(-\xi)^{-1}C\psi(-\xi A)B$, entailing~\eqref{eq.r-d}.
\end{pf}

Introducing the \emph{polylogarithm}~\cite{Wood92,wei2014} of order $s\in\mathbb{R}$ 
\begin{equation}\label{eq.polylog}
\mathrm{Li}_s(z)=\sum_{j=1}^{\infty}\frac{z^j}{j^s},\quad z\in\mathbb{C},|z|<1,
\end{equation}
it can be  checked that $\varphi(y)=-1+1/(1-e^{-y})=\mathrm{Li}_{0}(e^{-y})$ and, by using induction over $k$,
\begin{equation}\label{eq.phi-d-polylog}
\varphi^{(k)}(y) = (-1)^k\mathrm{Li}_{-k}(e^{-y}).
\end{equation}
\begin{rem}\label{rem.complete-mon} 
{\rm
Notice that $\mathrm{Li}_s(z)>0$ for $z$ being a real number from $(0,1)$.
Equality~\eqref{eq.phi-d-polylog} thus shows that $\varphi$ is \emph{completely monotonic}~\cite{MillerSamko2001}: $(-1)^k\varphi^{(k)}(y)>0$ for all $y>0$. In agreement with Corollary~\ref{cor.positivity1}, $r(\xi)>0$ for all $\xi>0$, thanks to~\eqref{eq.r} and~\eqref{eq.mean-value}, for any order $m$ and every choice of parameters $a_i,g_i>0$.}
\end{rem}

 As follows from Lemma~\ref{lem.Psi-props} below,  the function $\psi$, is \emph{not} completely monotonic, and hence~\eqref{eq.r-d}
 does not allow to establish that $r'(\xi)<0$ for all $\xi>0$. Nevertheless, for a low order $m$, the derivative $r'(\xi)$ is indeed sign-preserving, which allows to prove Theorem~\ref{thm.m9}. For an exact formulation, we state a corollary.

\begin{cor}\label{cor.monotone}
If $(-1)^{m-1}\psi^{(m-1)}(\zeta)>0$ at all $\zeta>0$, then $r'<0$ (i.e., $r$ is decreasing) on $(0,\infty)$. More generally,
$r$ is decreasing on any interval $(\xi_0,\xi_1)$ provided that $(-1)^{m-1}\psi^{(m-1)}(\zeta)>0$ for $\zeta\in (\xi_0\min_ia_i,\xi_1\max_ia_i)$.
\end{cor}
\begin{pf}
The proof is immediate from~\eqref{eq.r-d} and~\eqref{eq.mean-value} (applied to $k=m-1$) $\square$
\end{pf}
Corollary~\ref{cor.monotone} implies, e.g., that, for $m=3$, the function $r$ is decreasing, because $\psi$ is convex~\cite{Aut09}.
%In~\cite{Aut09}, the monotonicity of $r$ was derived from the Jenssen inequality, which approach cannot be easily generalized to $m>3$.

\subsubsection*{The derivatives $\psi^{(k)}$ and their estimates}

The derivatives of the function $\psi$, in fact, are also closely related to  polylogarithm~\eqref{eq.polylog} as summarized in the following lemma.

\begin{lem}\label{lem.Psi-props}
For each $k=1,2,\ldots$, one has
\begin{equation}\label{eq.Psi}
\begin{split}
\Psi_k(x) &\triangleq (-1)^k\psi^{(k)}(x)=\\&=x \mathrm{Li}_{-k-1}(e^{-x}) - k \mathrm{Li}_{-k}(e^{-x}),
%=\\=\sum_{j=1}^\infty j^k(x j - k)e^{-jx}.
\end{split}
\end{equation}
where $\Psi_k$  possesses the following properties:
\begin{enumerate}[{\rm (i)}]
\item $\Psi_k(x)>0$ for $0 \le x <\bar x(k)\triangleq 2 \pi/\sqrt[k+1]{2 k \boldsymbol{\zeta}(k+1)}$; 
\item $\Psi_k(x)>0$ for $x\geq k$;
\item $\Psi_k(x)>0$ for all $x>0$ if $k\leq\mcritprev$;
\item in general, the following inequality holds
\begin{equation}\label{eq.psi-bound}
\Psi_k(x)\geq -2k\frac{k!}{(2\pi)^{(k+1)}}\boldsymbol{\zeta}(k+1)\quad\forall x\geq 0,    
\end{equation}
where $\boldsymbol{\zeta}$ is the Riemann $\zeta$-function~\eqref{eq.zeta};
\item however, $\Psi_{\mcrit}$ attains negative values at some points $x>0$. 
\end{enumerate}
\end{lem}

\begin{figure}[h]
    \centering
    \includegraphics[width=0.45\textwidth]{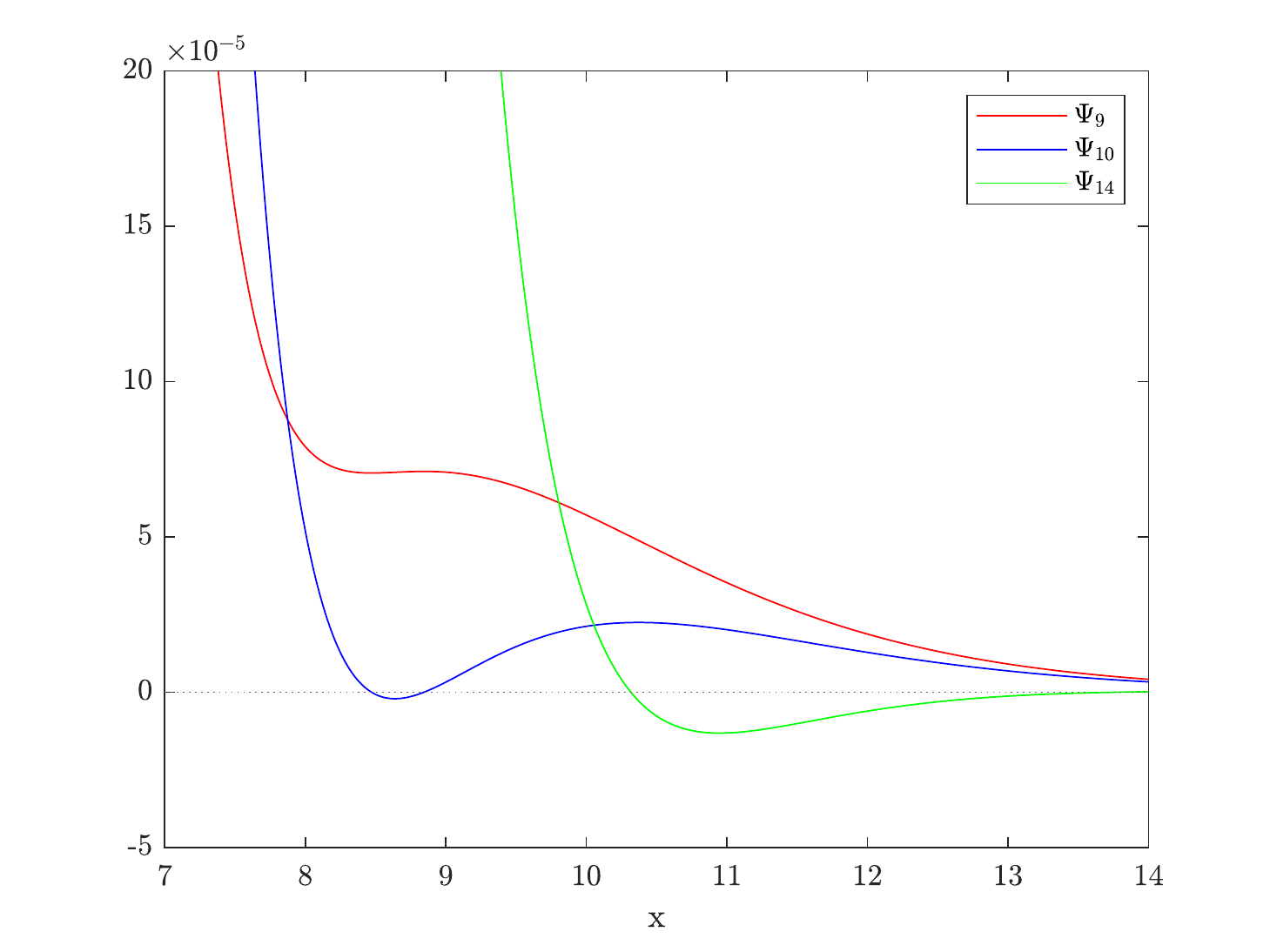}
    \caption{Graphs of $\Psi_k(x)$ for $k=9,10,14$.}
    \label{fig.Psi}
\end{figure}

The proof of Lemma~\ref{lem.Psi-props} is quite technical and given in Appendix~\ref{sec.app-psi}. 
Numerical simulation shows, in fact, that $\Psi_k$ is negative  for all $k>\mcritprev$. Fig.~\ref{fig.Psi} illustrates the behavior of $\Psi_{\mcritprev}$ (positive) and functions $\Psi_{\mcrit},\Psi_{14}$ that can attain negative values.

Combining the mean-value formula (Lemma~\ref{lem.mean.val}) with~\eqref{eq.r-d} and~\eqref{eq.psi-bound}, the following corollary is immediate.
\begin{cor}\label{cor.r-d-est}
For each $\xi>0$, the derivative $r'(\xi)$ admits the following upper bound
\begin{equation}\label{eq.aux1}
r'(\xi)<\frac{\bar g}{C_m},
\end{equation}
where $C_m$ is defined in~\eqref{eq.phi-d-bound}.
\end{cor}

As discussed in Section~\ref{sec.solution}, $C_m^{-1}$ decays exponentially as $m\to\infty$. Nevertheless, the right-hand side of~\eqref{eq.aux1} is positive, and 
for some parameters of IGO($m$), $m\geq 11$, it is possible that $r'$ attains positive values. 
In such a situation, the IGO may possess multiple 1-cycles (Section~\ref{sec.num}).

\subsubsection*{Proof of Theorem~\ref{thm.m9}}

If $m\leq\mcrit$, Lemma~\ref{lem.Psi-props} implies that $(-1)^{m-1}\psi^{(m-1)}(\xi)=\Psi_{m-1}(\xi)>0$ for all $\xi>0$;
in general, $\Psi_{m-1}(\xi)>0$ for $\xi\geq m-1$. Applying Corollary~\ref{cor.monotone}, one proves that $r(\cdot)$ decreases
on the interval $(\rho,\infty)$, where
\[
\rho\triangleq
\begin{cases}
0,\quad &1\leq m\leq 10,\\
\frac{m-1}{\min_ia_i},\,&\text{otherwise}.
\end{cases}
\]
When either $1\leq m\leq\mcrit$ or inequality~\eqref{eq.phi1} holds, then, obviously, $\Phi(y)>\rho$ for all $y\in\mathbb{R}_+$. 
Recalling that $\Phi,F$ are, respectively, non-decreasing and non-increasing, $R(y)=r(\Phi(y))F(y)$ is thus a non-increasing function, which means that equation~\eqref{eq.period} has only one solution on $\mathbb{R}_+$, and thus the IGO($m$) has a unique 1-cycle $\square$ 

\subsubsection*{Proof of Theorem~\ref{thm.m10}}

Theorem~\ref{thm.m10} is straightforward from Corollary~\ref{cor.r-d-est}. Indeed, the composition $r(\Phi(y))$ of a continuously differentiable (thus, locally Lipschitz) function and an absolutely continuous function is absolutely continuous, and one has
\[
\frac{d}{dy}r(\Phi(y))=r'(\Phi(y))\Phi'(y)<\frac{1}{F(0)},
\]
for almost all $y>0$ in view of~\eqref{eq.phi-d-bound} and~\eqref{eq.aux1}. The function $R(y)=r(\Phi(y))F(y)$ is now also absolutely continuous as a product of two absolutely continuous functions. Recalling that $F'(y)\leq 0$ at almost all $y>0$, one has $(y-R(y))'=1-r'(\Phi(y))\Phi'(y)F(y)-r(\Phi(y))F'(y)>0$, hence, $y-R(y)$ is increasing on $(0,\infty)$. 
Here, we used the fact that $0\leq F(y)\leq F(0)$ for all $y>0$
$\square$

%\subsubsection*{Additional remark}

\section{An example of the IGO with multiple 1-cycles.}\label{sec.num}

\label{The analytic result}

In this subsection, we construct IGO($m$) with at least three distinct 1-cycles for every $m$ such that
$\Psi_{m-1}(v_0)<0$ at some point $v_0$. This holds, e.g., for $m=\mcritnext$ (Lemma~\ref{lem.Psi-props}).

%We need the following auxiliary construction. Denote
Let $\Phi_{\sigma,y_*}$ be the Gaussian density distribution function with variance $\sigma^2$ and expectation $y_*$, that is,
\[
\Phi_{\sigma,y_*}(y)\triangleq\frac{1}{\sigma\sqrt{2\pi}}\int_{-\infty}^ye^{-\frac{(s-y_*)^2}{2\sigma^2}}ds.
\]
For each $\sigma>0$, one has $0<\Phi_{\sigma,y_*}(y_*)=1/2<\Phi_{\sigma,y_*}(\infty)=1$; also, $\Phi_{\sigma,y_*}$ is strictly increasing. By construction, the derivative 
\[
\Phi_{\sigma,y_*}'(y)=\frac{1}{\sigma\sqrt{2\pi}}e^{-\frac{(y-y_*)^2}{2\sigma^2}}
\]
attains its maximum $1/(\sigma\sqrt{2\pi})$ at $y=y_*$. 

The existence of multiple 1-cycles is established by the following lemma.
\begin{lem}\label{lem.counterex}
Choose numbers $y_*>0$, $\sigma>0$, and let $v_0>0$ be a point where $\Psi_{m-1}(v_0)<0$. Define IGO($m$) with the following parameters:
\begin{itemize}
\item a non-increasing differentiable function $F$ obeying~\eqref{eq:bounds};
\item $\Phi=\Phi_{\sigma,y_*}$ (this function is strictly increasing on $\mathbb{R}$); 
\item $a_1=\ldots=a_m=a\triangleq v_0/\Phi_{\sigma,y_*}(y_*)=2v_0$; 
\item finally, $g_1,\ldots,g_{m-1}>0$ are such that\footnote{Due to Remark~\ref{rem.complete-mon}, the right-hand side of~\eqref{eq.barg-aux} is positive.}
\begin{equation}\label{eq.barg-aux}
\bar g=(-1)^{(m-1)}\frac{2^{m-1}(m-1)!y_*}{\varphi^{(m-1)}(v_0)F(y_*)}. 
\end{equation}
\end{itemize}
Then, for a small enough $\sigma>0$, this IGO possesses at least three distinct positive 1-cycles.
\end{lem}
\begin{pf}
Combining~\eqref{eq.r} and~\eqref{eq.derivative}, one has 
\[
\begin{gathered}
r(\Phi(y_*))=\frac{(-1)^{(m-1)}\Phi(y_*)^{m-1}\bar g}{(m-1)!}\varphi^{(m-1)}(a(\sigma)\Phi(y_*))=\\=\frac{(-1)^{(m-1)}\bar g}{2^{m-1}(m-1)!}\varphi^{(m-1)}(v_0)=\frac{y_*}{F(y_*)}.
\end{gathered}
\]
Recalling that $R(y)=r(\Phi(y))F(y)$, one shows that $R(y_*)=y_*$.

Retracing the arguments from the proof of Theorem~\ref{thm.m10} above, one has
\[
\begin{split}
(y-R(y))'|_{y=y_*}=1-\underbrace{r'(\Phi(y_*))\Phi'(y_*)F(y_*)}_{P_1}-\\
-\underbrace{r(\Phi(y_*))F'(y_*)}_{P_2}.
\end{split}
\]
In view of~\eqref{eq.r} and~\eqref{eq.barg-aux}, $P_2$ does not depend on $\sigma$, being determined by $y_*$ and $F$ only:
\[
P_2=\frac{y_*F'(y_*)}{F(y_*)}.
\]
Recalling that $\Phi'(y_*)=1/(\sigma\sqrt{2\pi})$ and applying~\eqref{eq.r-d},
\[
\begin{gathered}
{P_1}=r'(\Phi(y_*))\Phi'(y_*)F(y_*)=\\
=\frac{(-1)^{m-2}\bar g\Phi(y_*)^{m-2}\psi^{(m-1)}(a(\sigma)\Phi(y_*))\Phi'(y_*)F(y_*)}{(m-1)!}=\\
=-\frac{\bar g F(y_*)\Psi_{m-1}(v_0)}{2^{m-2}(m-1)!\sigma\sqrt{2\pi}}=\\
=\frac{2(-1)^{m-1}y_*}{\varphi^{(m-1)}(v_0)}\frac{(-\Psi_{m-1}(v_0))}{\sigma\sqrt{2\pi}} >0.
\end{gathered}
\]
One notices that $P_1$ can be arbitrarily large for small $\sigma>0$; In particular, it is possible to choose $\sigma>0$ in such a way that $1-R'(y_*)<0$. Since $y_*-R(y_*)=0$, in there exists $\varepsilon\in (0,y_*)$  such that
\[
\begin{gathered}
y-R(y)>0,\quad y\in (y_*-\varepsilon,y_*),\\
y-R(y)<0,\quad y\in (y_*,y_*+\varepsilon).
\end{gathered}
\]
On the other hand, $y-R(y)<0$ as $y\to 0+$ and $y-R(y)\to+\infty$ as $y\to\infty$ (see the proof of Theorem~1). 
Hence,~\eqref{eq.period} has at least two additional solutions $y_1\in(0,y_*)$ and $y_2\in(y_*,\infty)$. 
In view of Theorem~\ref{thm.exist}, $y_1,y_*,y_2$ correspond to three distinct 1-cycles of the IGO
$\square$
\end{pf}

\begin{rem}
{\rm 
One may suspect that the existence of multiple 1-cycles is caused by the multiplicity of the eigenvalues $a_i=a=2v_0$, however, this is not the case.
The construct in Lemma~\ref{lem.counterex} can be generalized to the case where $a_i$ are close enough to $2v_0$ yet pairwise distinct. We omit this for brevity.
}
\end{rem}

\subsection{Numerical example}

The existence of multiple 1-cycles for the IGO of order $m=11$ is demonstrated now numerically by computations in Matlab, following the IGO construction method in Lemma~\ref{lem.counterex}. Set $y_*=2$, $\sigma=2\cdot 10^{-4}$, $F(y)=1$ (constant) and $v_0=8.64$, which corresponds to $\Psi_{m-1}(v_0)=\Psi_{10}(v_0)<0$.
We consider matrices~\eqref{eq:matrices}, where $a_1=\ldots=a_{11}=a=17.28$ and $g_1=\ldots=g_{10}=22.6486$, %That is correct /Håkan
which correspond, in view of~\eqref{eq.phi-d-polylog}, to the values
\[
\bar g=3.5515\cdot 10^{13}, \quad P_1=1.1257, \quad P_2=0.
\]
 In particular
\[
(y-R(y))'|_{y=y_*} = -0.1257,
\]
which indicates the existence of three solutions $y_1,y_*,y_2$ to the equation $y-R(y)=0$ and three corresponding 1-cycles. $y_1,y_2$ are found numerically to have the values
\[
y_1=1.9998234,\quad y_2=2.0002739.
\]
The fixed points of $Q(x)$ corresponding to $y_1,y_*,y_2$ are calculated according to Theorem~\ref{thm.exist} to
\[
x_1=\left[\begin{smallmatrix} 0.00019\\0.00216\\ 0.01213\\0.04535\\0.12731\\0.28635\\0.53835\\0.87280\\1.25282\\1.63477\\1.99982
\end{smallmatrix}\right],\; x_*=\left[\begin{smallmatrix} 0.00018\\0.00200\\0.01135\\0.04287\\0.12155\\0.27608\\0.52400\\0.85724\\1.24048\\1.62906\\2
\end{smallmatrix}\right],\; x_2=\left[\begin{smallmatrix} 0.00015\\0.00178\\0.01024\\0.03927\\0.11305\\0.26064\\0.50199\\0.83277\\1.22041\\1.61921\\2.00027
\end{smallmatrix}\right].
\]
Stability of the corresponding 1-cycles is determined by the Schur stability of the Jacobian matrix
\[
Q'(x) = e^{A \Phi(Cx)}(I+F'(Cx) BC) + \Phi'(Cx)A Q(x) C,
\]
evaluated at the fixed points. The numerical calculation shows that all three 1-cycles are unstable, and the spectral radii of
the corresponding Jacobian matrices are:
\[
\rho(Q'(x_1))=68.64,\rho(Q'(x_*))=64.91,\rho(Q'(x_2))=58.47. 
\]
%so all three 1-cycles are unstable.

\section{Conclusions}
 %         \textcolor{red}{TODO}    
 A special case of periodic solutions in the impulsive Goodwin's oscillator (IGO) characterized by one impulse generated by the pulse-modulated feedback in the least period, i.e. a 1-cycle, is considered. The continuous part of the IGO is allowed to be of arbitrary order, in contrast with the established in the literature case of third-order dynamics. The structure of the continuous part is still assumed to be a chain of first-order blocks. It is proved that a 1-cycle always exists in the IGO, regardless of the continuous part order. Further, when the continuous part order is at most ten, the 1-cycle is unique. It is demonstrated, by a constricting an example, that uniqueness does not generally apply to higher orders of the continuous part, e.g. for order eleven. Uniqueness of 1-cycle can however be recovered by restricting the slopes of the modulation functions of the IGO or even by restricting the feedback impulses to be sufficiently sparse. 
\bibliographystyle{plain}   
\bibliography{refs,observer}  
\appendix
%\section{Proof of the bounds~\eqref{eq.bounded1},\eqref{eq.bounded2}}\label{sec.app-bound}    % Each appendix must have a short title.
%\textcolor{red}{TODO if needed: can be removed (along with inequalities). The readers might believe in them.}

\section{Divided differences and Opitz formula}\label{sec.app-opitz}

Divided differences (DD) are widely used in numerical analysis and employed in this work to compute matrix functions. Here we review some basic properties of the DDs, referring the reader to~\cite{BerezinZhidkov,DeBoor2005,HornJohnsonTopics} for further details. 

\subsection*{Definitions of DD}

Throughout this section, we deal with functions $f:I\to\mathbb{R}$, where $I\subseteq\mathbb{R}$ is some interval (possibly, open).
The standard definition of the $k$-th order DD (briefly, $k$-DD) for such a function at a sequence of pairwise distinct points $x_0,\ldots,x_k\in I$ is as follows. We formally define the $0$-DD as $f[x_0]\triangleq f(x_0)$ and, subsequently, the $1$-DD as 
\[
f[x_0,x_1]\triangleq \frac{f(x_1)-f(x_0)}{x_1-x_0}.
\]
For $k\geq 2$, the $k$-DD is constructed inductively as
\begin{equation}\label{eq.k-dd}
f[x_0,\ldots,x_k]=\frac{f[x_1,\ldots,x_k]-f[x_0,\ldots,x_{k-1}]}{x_k-x_0}.
\end{equation}

An equivalent and more compact definition of the $k$-DD is based on the concept of \emph{interpolation polynomial}, which can be written in the Lagrange or Newton form. By definition, the interpolation polynomial of $f$ at the points $x_0,\ldots,x_k$ (where $x_i\ne x_j\,\forall i\ne j$) is the (unique) polynomial $L=L_{f,x_0,\ldots,x_k}$ of degree $\leq k$ such that
all $x_i$ are roots of the equation
\begin{equation}\label{eq.interp}
L(x)=f(x).
\end{equation}
It can be proven~\cite{BerezinZhidkov} that $L$ admits the form
\begin{equation}\label{eq.newton}
\begin{aligned}
L_{f,x_0,\ldots,x_k}&(x)=f[x_0]+f[x_0,x_1](x-x_0)+\ldots\\&+f[x_0,x_1,\ldots,x_k]\prod\nolimits_{j=0}^{k-1}(x-x_j),
\end{aligned}
\end{equation}
known as Newton's form of the interpolation polynomial. 
This leads to an alternative definition of the $k$-DD $f[x_0,x_1,\ldots,x_k]$, which is the \emph{lead} (degree $k$) coefficient of the interpolation polynomial. 

If $f$ is  differentiable  $k$ times on $I$, then the latter approach allows to define the $k$-DD to an \emph{arbitrary} sequence $x_0,\ldots,x_k$. If some number $\xi$ occurs $s$ times in this sequence ($1\leq s\leq k$), then $\xi$ is a root of~\eqref{eq.interp} with multiplicity $s$:
$
L^{(p)}(\xi)=f^{(p)}(\xi)\quad p=0,\ldots,s-1.
$
Adopting such a convention, the interpolation polynomial remains uniquely determined~\cite{BerezinZhidkov}, and hence its lead coefficient $f[x_0,x_1,\ldots,x_k]$ is well defined.

\textbf{Example:} If $x_0=\ldots=x_k=\xi$, then the interpolation polynomial is nothing else than the Taylor sum
\begin{equation}
L(x)=\sum_{j=0}^{k}\frac{f^{(j)}(\xi)}{j!}(x-\xi)^j, \label{eq.Lx}
\end{equation}
whose lead coefficient is 
\begin{equation*}%\label{eq.derivative}
f[\xi,\ldots,\xi]=f^{(k)}(\xi)/k!.
\end{equation*}

\subsection*{Technical properties of DDs}

%The expression in \eqref{eq.derivative} lends itself to the following important extension~\cite[Section~8]{DeBoor2005}, which relates the $k$-DD to the $k$-th %order derivative.

%\textcolor{blue}{Where is $s$ in the right-hand side of $\ell(s;x_0,\ldots,x_k)$?}
%\textcolor{red}{should be clear from the last line. Otherwise, formula is too long.}
In the next subsections, we will use the following simple property of the DD. 
%that, to the best of our knowledge, are not readily available in the literature.
\begin{lem}(Scaling)\label{lem.scale}
Given a function $f:(a,\infty)\to\mathbb{R}$ and a number $\xi\ne 0$, denote $f_{\xi}(x)\triangleq f(\xi x)$. Then
\[
f_{\xi}[x_0,\ldots,x_k]=\xi^kf[\xi x_0,\ldots,\xi x_k].
\]
\end{lem}
\begin{pf}
Notice that if $L(x)=L_{f,x_0,\ldots,x_k}$ is the interpolation polynomial for $f$, then $L(\xi x)$ is the interpolation polynomial for $f_{\xi}$. Recalling that $f_{\xi}[x_0,\ldots,x_k]$ and $f[x_0,\ldots,x_k]$ are the lead coefficients of respectively $L(\xi x)$, $L(x)$, one obtains the desired relation $\square$
\end{pf}

%\begin{lem}(Derivative with respect to a parameter).\label{lem.par}
%Consider a function $f(x,\eta)$ of two variables $x\in I$ and $\eta\in J$, where $I,J$ are some intervals. Denote the $k$-DD of %function $f(\cdot,\eta)$ (for $\eta$ being fixed) by $f[x_0,\ldots,x_k|\eta]$ and assume that the derivative $g(x,\eta):=%\frac{\partial}{\partial\eta}f(x,\eta)$ exists for all $x\in I,\eta\in J$. 
%Suppose also at the derivatives $\partial^k f(x,\eta)/\partial x^k$ and $\partial^{k+1} f(x,\eta)/\partial\eta\partial x^k$
%exist and are continuous on $I\times J$. Then
%\[
%\frac{\partial}{\partial\eta} f[x_0,\ldots,x_k|\eta]=g[x_0,\ldots,x_k|\eta].
%\]
%for all $x_0,\ldots,x_k\in I$, $\eta\in J$.
%\end{lem}
%\begin{pf}
%The proof follows\footnote{Notice that for pairwise distinct $x_i$, an alternative proof via induction on $k$ is possible based %on~\eqref{eq.k-dd}. In this situation, the existence of $k$th order derivative in $x$ is not needed.} from the Genocchi-Hermite %formula~\eqref{eq.genocchi} and standard theorems on differentiating under the integral sign\textcolor{red}{TODO: ref}.
%$\square$
%\end{pf}

Finally, we notice that the DDs linearly depend on $f$, that is, for two functions $f_1,f_2$ defined on $(a,b)$ and two coefficients $\alpha_1,\alpha_2$, one has $(\alpha_1f_1+\alpha_2f_2)[x_0,\ldots,x_k]=\alpha_1f_1[x_0,\ldots,x_k]+\alpha_2f_2[x_0,\ldots,x_k]$.

\subsection*{Functions on matrices and the Opitz formula}

%For a scalar function $f:\mathbb{C}\to\mathbb{C}$ that is holomorphic (complex-analytic) in  vicinity of the spectrum of the matrix $A$, the matrix function $f(A)$ can be defined, \cite[Section~6.2]{HornJohnsonTopics}. 
Let $\mathcal{D}\subseteq\mathbb{C}$ be an open simply connected set containing
the eigenvalues $\lambda_j$ of the matrix $A$ and $f:\mathbb{C}\to\mathbb{C}$  be holomorphic on $\mathcal{D}$. Then, for any simple closed curve $\Gamma\subset\mathcal{D}$ that encircles all $\lambda_j$ in the counter-clockwise direction \cite[Section~6.2]{HornJohnsonTopics},
\begin{equation}\label{eq:contour}
f(A)\triangleq \frac{1}{2\pi\imath}\oint_{\Gamma}f(z)(zI-A)^{-1}dz.
\end{equation}
%\textcolor{blue}{This is one of many equivalent matrix function definitions.}

In particular, if $S$ is an invertible matrix, then $f(SAS^{-1})=Sf(A)S^{-1}$. Also, for every two functions $f,g$, the matrices 
$f(A)$ and $g(A)$ commute. 

Consider now the two-diagonal matrix below 
\[
\Lambda=\begin{bmatrix}
\lambda_1 &0  &0  &0 \\
1 & \lambda_2 &0 &0 \\
0 &  \ddots & \ddots &0\\
0 &0 & 1 & \lambda_m
\end{bmatrix}.
%\quad \Lambda_{ij}:=
%\begin{cases}
%\lambda_i,&i=j\\
%1, &i=j+1,\\
%0,&\text{otherwise},
%\end{cases}
\]
Assuming $f$ complex analytic in  vicinity of $\lambda_1,\ldots,\lambda_n$, the matrix
$f(\Lambda)$ admits an elegant representation, known as the \emph{Opitz formula\footnote{Usually, the Opitz formula is given for upper-triangular two-diagonal matrices, the case of lower triangular is straightforward by noticing that $f(\Lambda^{\top})=f(\Lambda^{\top})^{\top}$.}}~\cite{Eller1987}.
Namely, $f(\Lambda)$ is the lower-triangular matrix whose entries are
\begin{equation}\label{eq.opitz}
%\begin{gathered}
(f(\Lambda))_{ij}=
\begin{cases}
f[\lambda_i,\ldots,\lambda_j],&i\geq j,\\
0,& i<j.
\end{cases}
%\end{gathered}
\end{equation}
For instance, the left-bottom corner entry is the $(m-1)$-DD of function $f$, that is, $f(\Lambda)_{m1}=f[\lambda_1,\ldots,\lambda_m]$.

\section{Proof of Lemma~\ref{lem.Psi-props}}\label{sec.app-psi}

To obtain the expression for $\Psi_k$, note that $\psi(x)$ can be expressed as a series:
\[
\begin{split}
\psi(x) &= \frac{xe^{x}}{(e^{x}-1)^2} = -x \bigg(\frac{1}{1-e^{-x}}\bigg)'= \\ &= -x \Big(\sum_{j=0}^\infty e^{-jx} \Big)' = \sum_{j=0}^\infty xj e^{-jx} = \sum_{j=1}^\infty xj e^{-jx},
\end{split}
\]
which implies the expression for the $k$-th derivative
\[
\begin{split}
\psi^{(k)}(x) &= (-1)^k \sum_{j=1}^\infty j^k(x j - k)e^{-jx}=\\
&= (-1)^k (x \mathrm{Li}_{-k-1}(e^{-x}) - k \mathrm{Li}_{-k}(e^{-x})),
\end{split}
\]
resulting in~\eqref{eq.Psi}.

Statement~(ii) follows from \cite[Theorem~8]{wei2014}.

To prove statements~(i) and~(iv), we need a representation of the polylogarithm of order $(-k)<0$~\cite[9.553]{gradshteyn2014}
\[
\mathrm{Li}_{-k}(e^{-x}) = k! \sum_{l=-\infty}^\infty (2 \pi li + x)^{-k-1},
\]
which leads to an alternative representation of $\Psi_k$:
\[
\begin{split}
\Psi_k(x)=x \mathrm{Li}_{-k-1}(e^{-x}) - k \mathrm{Li}_{-k}(e^{-x})=\\
(k+1)! \sum_{l=-\infty}^\infty x(2 \pi li + x)^{-k-2} - k k! \sum_{l=-\infty}^\infty (2 \pi li + x)^{-k-1}=\\
= k! \sum_{l=-\infty}^\infty ((k+1) x (2 \pi li + x)^{-k-2} - k (2 \pi li + x)^{-k-1})=\\
= k! \sum_{l=-\infty}^\infty (2 \pi li + x)^{-k-2} ((k+1) x - k(2 \pi li + x))=\\
%= k! \sum_{l=-\infty}^\infty (2 \pi li + x)^{-k-2}(x- 2 \pi kli)=\\
%= k! \sum_{l=-\infty}^\infty \frac{x+2 \pi li-2 \pi li - 2 \pi kli}{(2 \pi li + x)^{k+2}}=\\
= k!\big(\frac{1}{x^{k+1}} + \underbrace{\sum_{\substack{l=-\infty\\ l \ne 0}}^\infty \frac{x - 2 \pi kli}{(x + 2 \pi li)^{k+2}}}_{=h_k(x)}\big).
\end{split}
\]

Notice that for each $x>0$, one has
\[
\left|
\frac{x - 2 \pi kli}{(x + 2 \pi li)^{k+2}}
\right|\leq\frac{1}{|x + 2 \pi li|^{k+1}}\left|\frac{x - 2 \pi kli}{x + 2 \pi li}\right|,
\]
where the multipliers  are, obviously, less,  than $(2\pi |l|)^{-k-1}$ and $k$, respectively.
Statement~(iv) and~\eqref{eq.psi-bound} are now straightforward from the following estimate:
\[
|h_k(x)|\leq 2\sum_{l=1}^\infty\left|\frac{x - 2 \pi kli}{(x + 2 \pi li)^{k+2}}\right|\leq \frac{2k}{(2\pi)^{k+1}}\sum_{l=1}^{\infty}\frac{1}{l^{k+1}}.
\]

To prove statement~(i), it suffices to notice that $x=\bar x(k)$ is the unique real positive solution to the equation
\[
\frac{1}{x^{k+1}}-\frac{2k}{(2\pi)^{k+1}}\sum_{l=1}^{\infty}\frac{1}{l^{k+1}}=0;
\]
Obviously, $\Psi_k(x)>0$ as $0<x<\bar x(k)$. 

Statement~(iii) is proved similarly, refining the estimate for the term $h_k(x)$. Notice that for $k\leq 4$, statement (iii) follows from statements (i) and (ii), because $\bar x(k)<k$. For $k=5,\ldots,8$, one can use a more precise estimate:
\[
\begin{aligned}
h_k(x)&=\underbrace{\sum_{\substack{l=-\infty\\ l \ne 0,\pm 1}}^\infty \frac{x - 2 \pi kli}{(x + 2 \pi li)^{k+2}}}_{\tilde h_k(x)}+\\
&+\frac{x - 2 \pi ki}{(x + 2 \pi i)^{k+2}} + \frac{x + 2 \pi ki}{(x - 2 \pi i)^{k+2}},
\end{aligned}
\]
where $\tilde h_k(x)$ is estimated similarly to $h_k(x)$, that is,
\[
|\tilde h_k(x)|\leq \frac{2k}{(2\pi)^{k+1}}\sum_{l=2}^{\infty}\frac{1}{l^{k+1}}=
\frac{2k(\boldsymbol{\zeta}(k+1)-1)}{(2\pi)^{k+1}}.
%\underbrace{\sum_{l=2}^{\infty}\frac{1}{l^{k+1}}}_{=\boldsymbol{\zeta}(k+1)-1}
\]
Therefore, one obtains the following estimate for $\Psi_k$:
\[
\begin{split}
\frac{\Psi_k(x)}{k!} \ge \frac{1}{x^{k+1}} - \frac{2k}{(2\pi)^{k+1}} (\zeta(k+1)-1)+\\
+ \frac{x - 2 \pi ki}{(x + 2 \pi i)^{k+2}} + \frac{x + 2 \pi ki}{(x - 2 \pi i)^{k+2}}=\\
=\frac{p_k(x)}{x^{k+1}(x^2+4\pi^2)^{k+2}},
\end{split}
\]
where $p_k(x)$ is a polynomial of degree $3k+5$ such that $p_k(0)=4\pi^2>0$. To prove that $\Psi_k(x)>0$ for $x>0$, in view of statement (ii), it suffices to check that $p_k$ has no real roots on $[0,k]$. This is indeed the case for $k=5,\ldots,8$, as reported in Table~\ref{tab.polroots} (the roots were found numerically using Matlab), however, for $k=9$ this condition is violated.

To prove statement (iii) for $k=9$, one needs an even more refined estimate of $h_k$ as follows:
\[
\begin{aligned}
h_k(x)&=\underbrace{\sum_{\substack{l=-\infty\\ l \ne 0,\pm 1,\pm 2}}^\infty \frac{x - 2 \pi kli}{(x + 2 \pi li)^{k+2}}}_{\bar h_k(x)}+\\
&+\sum_{l=1}^2\frac{x - 2 \pi kli}{(x + 2 \pi li)^{k+2}} + \sum_{l=1}^2\frac{x + 2 \pi kli}{(x - 2 \pi li)^{k+2}},
\end{aligned}
\]
where $\bar h_k(x)$ can be estimated similar to $h_k,\tilde h_k$:
\[
|\bar h_k(x)|\leq \frac{2k}{(2\pi)^{k+1}}\sum_{l=3}^{\infty}\frac{1}{l^{k+1}}=
\frac{2k(\boldsymbol{\zeta}(k+1)-1-2^{-k-1})}{(2\pi)^{k+1}}.
\]
This entails a more refined estimate for $\Psi_k$:
\[
\begin{split}
\frac{\Psi_k(x)}{k!} \ge \frac{1}{x^{k+1}} - \frac{2k}{(2\pi)^{k+1}} (\zeta(k+1)-1-\frac{1}{2^{k+1}})\\
+ \frac{x - 2 \pi ki}{(x + 2 \pi i)^{k+2}} + \frac{x + 2 \pi ki}{(x - 2 \pi i)^{k+2}}\\
+ \frac{x - 4 \pi ki}{(x + 4 \pi i)^{k+2}} + \frac{x + 4 \pi ki}{(x - 4 \pi i)^{k+2}}=\\
=\frac{q_k(x)}{x^{k+1}(4\pi^2+x^2)^{k+2}(16\pi^2+x^2)^{k+2}},
\end{split}
\]
where $q_k(x)$ is a polynomial of order $5k+9$ satisfying $q_k(0)=(64\pi^2)^{k+2}>0$. As shown in Table~\ref{tab.polroots}, the real roots of $q_{\mcritprev}(x)$ are located outside the interval  $(0,\mcritprev]$.

Finally, statement (v) can be validated by computing the polylogarithmic functions in Matlab:
\[
\Psi_{\mcrit}(8.64)\approx -2.087496\cdot 10^{-6}\square
\]

\begin{table}%[ht]
    \centering
    \caption{Real roots of polynomials $p_k$ (for $k=5,\ldots,8$) and $q_9$.}
    \label{tab.polroots}
    \begin{tabular}{c|c}
        $k$ & Real roots \\
        \hline
        %$1$ & $\pm 4.178\dots$\\
        %$2$ & $5.915\dots$\\
        %$3$ & $\pm7.453\dots$\\
        %$4$ & $8.680\dots$\\
        $5$ & $\pm9.563\dots$\\
        $6$ & $10.115\dots$ \\
        $7$ & $\pm10.369\dots$\\
        $8$ & $10.291\dots$\\
        $9$ & $\pm15.456\dots$\\
        \hline
    \end{tabular}
    
\end{table}

\end{document}